# Transversality Conditions for Stochastic Higher-Order Optimality: Continuous and Discrete Time Problems

Dapeng Cai [a,*] and Takashi Gyoshin Nitta [b]

[a] Institute for Advanced Research, Nagoya University, Furo-cho, Chikusa-ku, Nagoya, 464-8601, Japan; [b] Department of Mathematics, Faculty of Education, Mie University, Kurimamachiya 1577, Tsu, 514-8507, Japan

**Abstract:** Higher-order optimization problems naturally appear when investigating the effects of a patent with finite length, as in the pioneering work of Futagami and Iwaisako (2007). In this paper, we establish the Euler equations and transversality conditions necessary for analyzing such higher-order optimization problems. We develop our results for stochastic general reduced-form models and consider cases of both continuous and discrete time. We employ our results to establish the Euler equations and transversality conditions for the simplified household maximization problem in Futagami and Iwaisako (2007).

---

[*] Corresponding Author. Fax: +81-52-788-6141. E-mail address: cai@iar.nagoya-u.ac.jp.

1. **Introduction**

In a pioneering work, Futagami and Iwaisako (2007) examined the effects of patent policies and showed that it is finite rather than infinite patent length that maximizes social welfare on the balanced growth path.[1] Optimization problems involving higher-order difference equations (particularly the household maximization problem) naturally appear in their analysis. However, Futagami and Iwaisako (2007) did not explicitly consider the transversality conditions (hereafter, TVCs), even though both the Euler equations and TVCs would be necessary to identify the optimums in the optimization problems. In fact, to our best knowledge, readily applicable general results on the Euler equations and TVCs for higher-order optimization problems are not available in the literature. This is despite optimization problems involving higher-order differential or difference equations becoming increasingly important in dynamic economic analysis.[2] For example, recent discussion concerning intertemporally dependent preferences, in which the agent's utility in a particular period depends on both current as well as previous consumption levels, also involves higher-order optimization problems (Mitra and Nishimura 2006). The purpose of this paper is to develop conditions that ensure the necessity of the TVCs for higher-order optimization problems, even when the objective functional is stochastic and unbounded.

---

[1] The model Futagami and Iwaisako (2007) considered does not exhibit scale effects.

[2] Economic models also include higher-order differential problems. More specially, second-order differential problems appear in discussions on the acceleration principle, in that when investment depends on the variation in income, consumption depends on the second-order differential of capital.

In this paper, we use the elementary variational approach to derive the Euler equations and the TVCs for stochastic infinite horizon optimality, and consider both continuous and discrete time problems. This approach enables us to obtain the TVCs for extremely general reduced-form models with very few technical restrictions. Other than the standard assumptions, we identify those assumptions necessary when the objective functionals do not converge.[3] We also consider two simple examples. In particular, we show that depending on the properties of the perturbation curves under consideration, the examples can also represent counterexamples, in which the assumptions are violated and the TVCs are not satisfied. Finally, we apply our results to establish the Euler equations and the TVCs for the simplified household maximization problem in Futagami and Iwaisako (2007).

Our results should be useful for economic analyses of models with unbounded returns that involve higher-order differential or difference equations. They thus apply to, say, formal economic analyses of sustainability and global environmental problems, as these naturally involve higher-order difference equations, unbounded returns, and uncertainty. This is because of the following reasons. First, higher-order difference equations appear because mitigation efforts and investments generally take several decades to complete, and there can be significant time lags before any effects are noted. Second, discussion of intergenerational equity, especially when concerning the choice of discount factors, may cause objective functionals to diverge. Finally, climate policies must be undertaken when the impact of global warming is not fully known, and when the future economic and social consequences of climate change, particularly the valuation of possible damage, are uncertain.

---

[3] For approaches on how to explicitly construct optimal solutions to infinite horizon optimization problems with unbounded objective functionals, see Cai and Nitta (2009, 2011).

There is a growing literature on the necessity of TVCs. In one important contribution, Kamihigashi (2001) employs the Lebesgue integral to present the most general form of TVCs for first-order differential problems, $\max_{\mathbf{x}} \int_0^\infty v(\mathbf{x}(t), \dot{\mathbf{x}}(t), t) dt$. This generalizes the results in Weitzman (1973), Ekeland and Scheinkman (1986) and Michel (1990).[4] In other work using the classical Riemannian integral, Okumura, Cai and Nitta (2009) further extend these results to higher-order differential problems. Cai and Nitta (2010) also present the TVCs for deterministic higher-order discrete time problems. Alternatively, Kamihigashi (2003) represents an extension of Kamihigashi (2001) to the discrete time stochastic case.[5] In this paper, we extend the extant results to stochastic higher-order problems. As in Kamihigashi (2001, 2003), we consider extremely general reduced-form models.

The remainder of the paper is as follows. Section 2 presents the model and states the main results. In Section 3, we discuss two examples using continuous and discrete time. We show that these examples can also become counterexamples depending on the properties of the perturbation curves under consideration. Section 4 investigates the correspondence between the results for the continuous and discrete time models. Section 5 applies our results to the simplified household maximization problem examined in Futagami and Iwaisako (2007). Finally, Section 6 concludes the paper by discussing, among other things, how our results relate to previous work in the area.

---

[4] See Kamihigashi (2000, 2001, 2002, 2003) for the discussion and extension of these results.

[5] As argued in Kamihigashi (2003), the stochastic versions of Weitzman's (1973) theorem by Zilcha (1976) and Takekuma (1992) are not easily applicable as they use support prices and thus depend heavily on the infinite-dimensional separation theorem, which imposes several severe restrictions.

## 2. Main Results

### 2.1 *Derivation of the TVCs for the Continuous Time Problems*

We first use the elementary variational approach to present a complete characterization of both the Euler equations and TVCs for the stochastic higher-order continuous time optimization problems.

Let $(\Omega, F, P)$ be a probability space. Let $E$ denote the associated expectation operator, i.e., $Ez = \int z(\omega)dP(\omega)$ for any random variable $z: \Omega \to \bar{\mathbb{R}}$. We consider the following problem:

$$\begin{cases} \max_{\mathbf{x}} \int_0^\infty Ev\big(\mathbf{x}(t,\omega), \dot{\mathbf{x}}(t,\omega), \ddot{\mathbf{x}}(t,\omega), \cdots, \mathbf{x}^{(n)}(t,\omega), t, \omega\big) dt \\ \text{subject to } \mathbf{x}(0,\omega) = \mathbf{x}_0(\omega), \\ \forall t \geq 0, \; \big(\mathbf{x}(t,\omega), \dot{\mathbf{x}}(t,\omega), \ddot{\mathbf{x}}(t,\omega), \cdots, \mathbf{x}^{(n)}(t,\omega)\big) \in X(t,\omega) \subset \big(\mathbb{R}^N\big)^{n+1}, \end{cases} \quad (1)$$

where $N \in \mathbb{N}$, $v$ is a real-valued $n$th-order continuously differentiable function, and $\mathbf{x}(t,\omega)$ is $n$th-order continuously differentiable.[6] Notice that the objective functional of (1) is not necessarily finite. Let $v_1(y,z)$ denote the partial derivative of $v$ with respect to $y$; define $v_2(y,z)$ similarly.

Following Kamihigashi (2003), we make several standard assumptions. We assume that there exists a sequence of real vector space $\{B_t\}_{t=0}^{\infty}$ such that $\bar{\mathbf{x}}_0 \in F(\Omega, B_0)$ and $\forall t \in \mathbb{R}_+$, $X_t \subset F(\Omega, B_t) \times F(\Omega, B_{t+1}) \times \cdots \times F(\Omega, B_{t+N-1})$. Moreover, $\forall t \in \mathbb{R}_+$, $\forall (\mathbf{x}(t,\omega), \dot{\mathbf{x}}(t,\omega), \ddot{\mathbf{x}}(t,\omega), \cdots, \mathbf{x}^{(n)}(t,\omega)) \in X(t,\omega)$, we assume that

(i) $\forall \omega \in \Omega$, $v(\mathbf{x}(t,\omega), \dot{\mathbf{x}}(t,\omega), \ddot{\mathbf{x}}(t,\omega), \cdots, \mathbf{x}^{(n)}(t,\omega), t, \omega) \in [-\infty, \infty)$,

(ii) the mapping $v(\mathbf{x}(t,\omega), \dot{\mathbf{x}}(t,\omega), \ddot{\mathbf{x}}(t,\omega), \cdots, \mathbf{x}^{(n)}(t,\omega), t, \omega) : \Omega \to [-\infty, \infty)$ is measurable, and

(iii) $Ev(\mathbf{x}(t,\omega), \dot{\mathbf{x}}(t,\omega), \ddot{\mathbf{x}}(t,\omega), \cdots, \mathbf{x}^{(n)}(t,\omega), t, \omega)$ exists in $[-\infty, \infty)$.

We also assume that the optimal path to (1) exists and is given by $\mathbf{x}^*(t,\omega)$, optimal in the sense of an overtaking criterion.[7] We perturb it with $n$th-order continuously differentiable curves $\mathbf{p}(t,\omega)$,

$$\mathbf{x}(t,\omega) = \mathbf{x}^*(t,\omega) + \varepsilon \cdot \mathbf{p}(t,\omega). \tag{2}$$

We define

$$V(\varepsilon, T) = \inf_{T \leq T'} \int_0^{T'} (Ev(\mathbf{x}^*(t,\omega) + \varepsilon \cdot \mathbf{p}(t,\omega), \dot{\mathbf{x}}^*(t,\omega) + \varepsilon \dot{\mathbf{p}}(t,\omega), \cdots, \mathbf{x}^{*(n)}(t,\omega) + \varepsilon \cdot \mathbf{p}^{(n)}(t,\omega), t, \omega)$$
$$- Ev(\mathbf{x}^*(t,\omega), \dot{\mathbf{x}}^*(t,\omega), \cdots, \mathbf{x}^{*(n)}(t,\omega), t, \omega))dt.$$

---

[6] Normally, $v$ is defined on $(\mathbb{R}^N)^{n+1} \times \mathbb{R} \times \Omega$. The domain of $v$ is denoted by $X(t,\omega)$, included in $(\mathbb{R}^N)^{n+1}$, for all $t$, $\omega$.

[7] Clearly, our analysis only serves for the case where the optimal solution is interior to $X(t,\omega)$.

(3)

We assume that the optimal path satisfies the weak maximality criterion, á la Brock (1970):[8]

**Definition 1.** An attainable path $\{\mathbf{x}^*(t,\omega)\}$ is optimal if no other attainable path overtakes it:

$$\liminf_{T \to \infty \; T \leq T'} \int_0^{T'} (Ev(\mathbf{x}^*(t,\omega) + \varepsilon \cdot \mathbf{p}(t,\omega), \dot{\mathbf{x}}^*(t,\omega) + \varepsilon \dot{\mathbf{p}}(t,\omega), \cdots, \mathbf{x}^{*(n)}(t,\omega) + \varepsilon \cdot \mathbf{p}^{(n)}(t,\omega), t, \omega) \quad (4)$$
$$- Ev(\mathbf{x}^*(t,\omega), \dot{\mathbf{x}}^*(t,\omega), \cdots, \mathbf{x}^{*(n)}(t,\omega), t, \omega))dt \leq 0.$$

Let $V(\varepsilon) = \lim_{T \to \infty} V(\varepsilon, T)$. Differentiating it with respect to $\varepsilon$, we have

$$\lim_{\varepsilon \to {}^+0} \frac{V(\varepsilon)}{\varepsilon}$$
$$= \lim_{\varepsilon \to {}^+0} \liminf_{T \to \infty \; T \leq T'} \int_0^{T'} \frac{1}{\varepsilon} (Ev(\mathbf{x}^*(t,\omega) + \varepsilon \cdot \mathbf{p}(t,\omega), \dot{\mathbf{x}}^*(t,\omega) + \varepsilon \dot{\mathbf{p}}(t,\omega), \cdots, \mathbf{x}^{*(n)}(t,\omega) + \varepsilon \cdot \mathbf{p}^{(n)}(t,\omega), t, \omega)$$
$$- Ev(\mathbf{x}^*(t,\omega), \dot{\mathbf{x}}^*(t,\omega), \cdots, \mathbf{x}^{*(n)}(t,\omega), t, \omega))dt.$$

(5)

Let $\lim_{\varepsilon \to {}^+0} \frac{V(\varepsilon)}{\varepsilon} \equiv \Pi$. In general, $\frac{d}{d\varepsilon} \lim_{T \to \infty} f(\varepsilon, T) = \lim_{T \to \infty} \frac{d}{d\varepsilon} f(\varepsilon, T)$ only when $\lim_{T \to \infty} \frac{d}{d\varepsilon} f(\varepsilon, T)$ converges uniformly for $\varepsilon$ (Lang, 1997). As in Okumura, Cai and Nitta (2009), we impose the following two assumptions, which correspond to Assumption 3.1 in Kamihigashi (2001) when $n = 1$, although instead of uniform convergence, Kamihigashi (2001) used the definition of sup norm.

---

[8] Brock (1970) showed that once such a path exists once two assumptions are satisfied.

**Assumption 1.** *Assume $\Pi$ converges uniformly for $\varepsilon$ when $T \to \infty$.*

**Assumption 2.** *We assume that for any $T > 0$,*

$$\inf_{T \leq T'} \int_0^{T'} \frac{\left(Ev\left(\mathbf{x}^*(t,\omega) + \varepsilon \cdot \mathbf{p}(t,\omega), \dot{\mathbf{x}}^*(t,\omega) + \varepsilon \cdot \dot{\mathbf{p}}(t,\omega), \cdots, t, \omega\right) - Ev\left(\mathbf{x}^*(t,\omega), \dot{\mathbf{x}}^*(t,\omega), \cdots, t, \omega\right)\right)}{\varepsilon} dt$$

*converges uniformly for $\varepsilon$.*

Let

$$A(T', \varepsilon) = \int_0^{T'} \frac{\left(Ev\left(\mathbf{x}^*(t,\omega) + \varepsilon \cdot \mathbf{p}(t,\omega), \dot{\mathbf{x}}^*(t,\omega) + \varepsilon \cdot \dot{\mathbf{p}}(t,\omega), \cdots, t, \omega\right) - Ev\left(\mathbf{x}^*(t,\omega), \dot{\mathbf{x}}^*(t,\omega), \cdots, t, \omega\right)\right)}{\varepsilon} dt.$$

Assumption 2 means that there exists a sequence $A(T'_n, \varepsilon)$ for each $\varepsilon > 0$ such that $\lim_{n \to \infty} A(T'_n, \varepsilon) = \inf_{T \leq T'} A(T', \varepsilon)$, uniformly for $\varepsilon$, that is, the sequence is uniformly convergent for $\varepsilon$ (Okumura, Cai and Nitta, 2009). If $\Pi$ satisfies Assumption 1, we can then rewrite (5) as

$$\Pi = \lim_{T \to \infty} \lim_{\varepsilon \to {}^+0} \inf_{T \leq T'} A(T', \varepsilon). \tag{6}$$

On the other hand, when Assumption 2 is satisfied, $\lim_{\varepsilon \to {}^+0}$ and $\inf_{T \leq T'}$ can be interchanged, and equality (6) can be further restated as

$$\Pi = \lim_{T \to \infty} \inf_{T \leq T'} \lim_{\varepsilon \to {}^+0} A(T', \varepsilon). \tag{7}$$

Because $T'$ is finite uniformly for $\varepsilon$, if $A(T', \varepsilon)$ exists, we can then restate (7) as

$$\Pi = \lim_{T \to \infty} \inf_{T \leq T'} \int_0^{T'} \lim_{\varepsilon \to {}^+0} \frac{1}{\varepsilon}(Ev(\mathbf{x}^*(t,\omega) + \varepsilon \cdot \mathbf{p}(t,\omega), \dot{\mathbf{x}}^*(t,\omega) + \varepsilon \dot{\mathbf{p}}(t,\omega), \cdots, \mathbf{x}^{(n)*}(t,\omega) + \varepsilon \cdot \mathbf{p}^{(n)}(t,\omega), t, \omega)$$
$$- Ev(\mathbf{x}^*(t,\omega), \dot{\mathbf{x}}^*(t,\omega), \cdots, \mathbf{x}^{(n)*}(t,\omega), t, \omega))dt.$$

(7')

Since $v$ is differentiable, we see that

$$\lim_{\varepsilon \to^+ 0} \frac{1}{\varepsilon}(Ev(\mathbf{x}^*(t,\omega)+\varepsilon\cdot\mathbf{p}(t,\omega),\dot{\mathbf{x}}^*(t,\omega)+\varepsilon\dot{\mathbf{p}}(t,\omega),\cdots,\mathbf{x}^{(n)*}(t,\omega)+\varepsilon\cdot\mathbf{p}^{(n)}(t,\omega),t,\omega)$$
$$-Ev(\mathbf{x}^*(t,\omega),\dot{\mathbf{x}}^*(t,\omega),\cdots,\mathbf{x}^{(n)*}(t,\omega),t,\omega)) \tag{8}$$
$$=\frac{d}{d\varepsilon}(Ev(\mathbf{x}^*(t,\omega)+\varepsilon\cdot\mathbf{p}(t,\omega),\dot{\mathbf{x}}^*(t,\omega)+\varepsilon\dot{\mathbf{p}}(t,\omega),\cdots,\mathbf{x}^{(n)*}(t,\omega)+\varepsilon\cdot\mathbf{p}^{(n)}(t,\omega),t,\omega)$$
$$-Ev(\mathbf{x}^*(t,\omega),\dot{\mathbf{x}}^*(t,\omega),\cdots,\mathbf{x}^{(n)*}(t,\omega),t,\omega)).$$

We impose another assumption:

**Assumption 3.** Define

$$m_t(\varepsilon,\omega) \equiv \frac{1}{\varepsilon}(v(\mathbf{x}^*(t,\omega)+\varepsilon\cdot\mathbf{p}(t,\omega),\dot{\mathbf{x}}^*(t,\omega)+\varepsilon\dot{\mathbf{p}}(t,\omega),\cdots,\mathbf{x}^{(n)*}(t,\omega)+\varepsilon\cdot\mathbf{p}^{(n)}(t,\omega),t,\omega)$$
$$-v(\mathbf{x}^*(t,\omega),\dot{\mathbf{x}}^*(t,\omega),\cdots,\mathbf{x}^{(n)*}(t,\omega),t,\omega)).$$

Assume that there exists $\hat{m}_t(\tilde{\varepsilon},\omega)$, such that $|m_t(\varepsilon,\omega)| < \hat{m}_t(\tilde{\varepsilon},\omega)$, $\varepsilon \in (0,\tilde{\varepsilon}]$.

Clearly, Assumption 3 corresponds to Assumption (A.4) in Kamihigashi (2003) when $n=1$.

**Theorem 1.** Under Assumption 1-3, for any interior optimal path $\mathbf{x}^*(t,\omega)$, the Euler equation is given by

$$v_1\left(\mathbf{x}^*(t,\omega),\dot{\mathbf{x}}^*(t,\omega),\cdots,\mathbf{x}^{*(n)}(t,\omega),t,\omega\right) - \left(v_2\left(\mathbf{x}^*(t,\omega),\dot{\mathbf{x}}^*(t,\omega),\cdots,\mathbf{x}^{*(n)}(t,\omega),t,\omega\right)\right)' $$
$$+\cdots+(-1)^n\left(v_{n+1}\left(\mathbf{x}^*(t,\omega),\dot{\mathbf{x}}^*(t,\omega),\cdots,\mathbf{x}^{*(n)}(t,\omega),t\right)\right)^{(n)} = 0, \tag{9}$$

whereas the TVC is given by

$$\liminf_{T\to\infty} E\left[p\left(v_2-(v_3)'+\cdots+(-1)^{n-1}(v_n)^{(n-1)}\right)+\dot{p}\left(v_3-(v_4)'+\cdots+(-1)^{n-2}(v_{n-1})^{(n-2)}\right)+\cdots+p^{(n-1)}v_{n+1}\right]_0^T \leq 0.$$

$$\tag{10}$$

**Proof.** The right hand side (RHS) of (8) can be stated as

$$\frac{d}{d\varepsilon}(Ev(\mathbf{x}^*(t,\omega)+\varepsilon\cdot\mathbf{p}(t,\omega),\dot{\mathbf{x}}^*(t,\omega)+\varepsilon\dot{\mathbf{p}}(t,\omega),\cdots,\mathbf{x}^{(n)*}(t,\omega)+\varepsilon\cdot\mathbf{p}^{(n)}(t,\omega),t,\omega)$$
$$-Ev(\mathbf{x}^*(t,\omega),\dot{\mathbf{x}}^*(t,\omega),\cdots,\mathbf{x}^{(n)*}(t,\omega),t,\omega))$$
$$=\lim_{\varepsilon^+\to\infty}\frac{1}{\varepsilon}(Ev(\mathbf{x}^*(t,\omega)+\varepsilon\cdot\mathbf{p}(t,\omega),\dot{\mathbf{x}}^*(t,\omega)+\varepsilon\dot{\mathbf{p}}(t,\omega),\cdots,\mathbf{x}^{(n)*}(t,\omega)+\varepsilon\cdot\mathbf{p}^{(n)}(t,\omega),t,\omega)$$
$$-Ev(\mathbf{x}^*(t,\omega),\dot{\mathbf{x}}^*(t,\omega),\cdots,\mathbf{x}^{(n)*}(t,\omega),t,\omega))$$
$$=\lim_{\varepsilon^+\to\infty}\frac{1}{\varepsilon}(\int v(\mathbf{x}^*(t,\omega)+\varepsilon\cdot\mathbf{p}(t,\omega),\dot{\mathbf{x}}^*(t,\omega)+\varepsilon\dot{\mathbf{p}}(t,\omega),\cdots,\mathbf{x}^{(n)*}(t,\omega)+\varepsilon\cdot\mathbf{p}^{(n)}(t,\omega),t,\omega)dp(\omega)$$
$$-\int v(\mathbf{x}^*(t,\omega),\dot{\mathbf{x}}^*(t,\omega),\cdots,\mathbf{x}^{(n)*}(t,\omega),t,\omega))dp(\omega))$$
$$=\lim_{\varepsilon^+\to\infty}(\int\frac{1}{\varepsilon}(v(\mathbf{x}^*(t,\omega)+\varepsilon\cdot\mathbf{p}(t,\omega),\dot{\mathbf{x}}^*(t,\omega)+\varepsilon\dot{\mathbf{p}}(t,\omega),\cdots,\mathbf{x}^{(n)*}(t,\omega)+\varepsilon\cdot\mathbf{p}^{(n)}(t,\omega),t,\omega)$$
$$-v(\mathbf{x}^*(t,\omega),\dot{\mathbf{x}}^*(t,\omega),\cdots,\mathbf{x}^{(n)*}(t,\omega),t,\omega))dp(\omega)).$$

When Assumption 3 is satisfied, following Lebesgue's dominated convergence theorem, we then have

$$\lim_{\varepsilon^*\to\infty}(\int\frac{1}{\varepsilon}(v(\mathbf{x}^*(t,\omega)+\varepsilon\cdot\mathbf{p}(t,\omega),\dot{\mathbf{x}}^*(t,\omega)+\varepsilon\dot{\mathbf{p}}(t,\omega),\cdots,\mathbf{x}^{(n)*}(t,\omega)+\varepsilon\cdot\mathbf{p}^{(n)}(t,\omega),t,\omega)$$
$$-v(\mathbf{x}^*(t,\omega),\dot{\mathbf{x}}^*(t,\omega),\cdots,\mathbf{x}^{(n)*}(t,\omega),t,\omega))dp(\omega))$$
$$=E\frac{d}{d\varepsilon}(v(\mathbf{x}^*(t,\omega)+\varepsilon\cdot\mathbf{p}(t,\omega),\dot{\mathbf{x}}^*(t,\omega)+\varepsilon\dot{\mathbf{p}}(t,\omega),\cdots,\mathbf{x}^{(n)*}(t,\omega)+\varepsilon\cdot\mathbf{p}^{(n)}(t,\omega),t,\omega)$$
$$-v(\mathbf{x}^*(t,\omega),\dot{\mathbf{x}}^*(t,\omega),\cdots,\mathbf{x}^{(n)*}(t,\omega),t,\omega))$$
$$=E(v_1\mathbf{p}(t,\omega)+v_2\dot{\mathbf{p}}(t,\omega)+\cdots+v_{n+1}\mathbf{p}^{(n)}(t,\omega)).$$

Hence, $\Pi = \liminf\limits_{T\to\infty \atop T\le T'}\int_0^{T'} E\left(v_1\mathbf{p}(t,\omega)+v_2\dot{\mathbf{p}}(t,\omega)+\cdots+v_{n+1}\mathbf{p}^{(n)}(t,\omega)\right)dt.$

Using partial integration, we obtain

$$\int_0^T E\left(v_{1^{(k)}}p^{(k)}\right)dt = E\left(p^{(k-1)}v_{2^{(k)}}\right) - \int_0^T E\left(p^{(k-1)}(v_{2^{(k)}})'\right)dt, \quad (11)$$

and

$$\int_0^T E(v_1\mathbf{p}(t,\omega)+v_2\dot{\mathbf{p}}(t,\omega)+\cdots+v_{n+1}\mathbf{p}^{(n)}(t,\omega))dt$$
$$=\int_0^T E(v_1-(v_2)'+\cdots+(-1)^n(v_{n+1})^{(n)})p(t,\omega)dt + E[p(t,\omega)(v_2-(v_3)'+\cdots+(-1)^{n-1}(v_{n+1})^{(n-1)})$$
$$+\dot{p}(t,\omega)(v_3-(v_4)'+\cdots+(-1)^{n-2}(v_{n+1})^{(n-2)})+\cdots+p^{(n-1)}(t,\omega)v_{n+1}]_0^T.$$

Note that the argument would be the same when $\varepsilon \to {}^-0$. Therefore, for an arbitrary $n$th-order continuously differentiable curve $\mathbf{p}(t,\omega)$ satisfying Assumptions 1-3 (for example, when $\mathbf{p}(t,\omega)$ has compact support), we then have

$$\begin{aligned} 0 &\geq \liminf_{T\to\infty} \int_0^T \frac{d}{d\varepsilon} Ev(t, \mathbf{x}+\varepsilon p, \cdots, \mathbf{x}^{(n)} + \varepsilon p^{(n)}) dt \\ &\geq \liminf_{T\to\infty} \int_0^T E(v_1 - (v_2)' + \cdots + (-1)^n (v_{n+1})^{(n)}) p(t,\omega) dt \\ &\quad + \liminf_{T\to\infty} E[p(t,\omega)(v_2 - (v_3)' + \cdots + (-1)^{n-1}(v_{n+1})^{(n-1)}) + \dot{p}(t,\omega)(v_3 - (v_4)' \\ &\quad + \cdots + (-1)^{n-2}(v_{n+1})^{(n-2)}) + \cdots + p^{(n-1)}(t,\omega) v_{n+1}]_0^T. \end{aligned} \tag{12}$$

*Q.E.D.*

Clearly, (9) generalizes the standard Euler equation $v_1 - (v_2)' = 0$, whereas (10) is a stochastic version of the TVC in Okumura, Cai and Nitta (2009). Next we consider the linkage between our results and that of Kamihigashi (2001). We fix $0 < \bar{\alpha} < 1$ and $\alpha : \mathbb{R}^+ \to \mathbb{R}^+$, $C^{(n-1)}$, $\alpha(0) = 0, \cdots, \alpha^{(n-1)}(0) = 0$, $\alpha(t) = \bar{\alpha}$, $t \geq 1$. We consider a special curve $\mathbf{p}(t,\omega)$. Let $p(t,\omega) = \alpha x^*(t,\omega)$, then (10) is modified to

$$\begin{aligned} &\liminf_{T\to\infty} E[\alpha x^*(v_2 - (v_3)' + \cdots + (-1)^{n-1}(v_{n+1})^{(n-1)}) + (\alpha x^*)'(v_2 - (v_3)' \\ &\quad + \cdots + (-1)^{n-2}(v_{n+1})^{(n-2)}) + \cdots + (\alpha x^*)^{(n-1)} v_{n+1}]_0^T \\ &= \bar{\alpha} \liminf_{T\to\infty} E(x^*(v_2 - (v_3)' + \cdots + (-1)^{n-1}(v_{n+1})^{(n-1)}) + (x^*)'(v_2 - (v_3)' \\ &\quad + \cdots + (-1)^{n-2}(v_{n+1})^{(n-2)}) + \cdots + (x^*)^{(n-1)} v_{n+1}) \Big|_T \\ &\leq 0. \end{aligned}$$

Because $\bar{\alpha} > 0$, we have

$$\begin{aligned} &\liminf_{T\to\infty} E(x^*(v_2 - (v_3)' + \cdots + (-1)^{n-1}(v_{n+1})^{(n-1)}) + (x^*)'(v_3 - (v_4)' + \cdots \\ &\quad + (-1)^{n-2}(v_{n+1})^{(n-2)}) + \cdots + (x^*)^{(n-1)} v_{n+1}) \Big|_T \leq 0, \end{aligned}$$

which is then a stochastic extension of Kamihigashi (2001)'s TVC.

## 2.2  Derivation of the TVCs for the Discrete Time Problems

We proceed to consider the following stochastic higher-order difference problem:

$$\begin{cases} \max_{\mathbf{y}} \sum_{t=0}^{\infty} EV\left(\mathbf{y}(t,\omega),\mathbf{y}(t+1,\omega),\cdots,\mathbf{y}(t+n,\omega),t,\omega\right) \\ \text{subject to } \mathbf{y}(0,\omega) = \mathbf{y}_0(\omega), \ \forall t \geq 0, \\ \left(\mathbf{y}(t,\omega),\mathbf{y}(t+1,\omega),\cdots,\mathbf{y}(t+N-1,\omega)\right) \in X(t,\omega) \subset \left(\mathbb{R}^N\right)^{n+1}, \end{cases} \quad (13)$$

where $N \in \mathbb{N}$, and $V$ is a real-valued first-order continuously differentiable function.[9] Notice that the objective functional of (13) can be infinite.

Again, as in Kamihigashi (2003, Assumption 2.1 and 2.2), we assume that there exists a sequence of real vector space $\{B_t\}_{t=0}^{\infty}$ such that $\overline{\mathbf{y}}_0 \in F(\Omega, B_0)$ and $\forall t \in \mathbb{Z}_+, X(t) \subset F(\Omega, B_t) \times F(\Omega, B_{t+1}) \times \cdots \times F(\Omega, B_{t+n})$. Moreover, $\forall t \in \mathbb{Z}_+$, $\forall \left(\mathbf{y}(t,\omega),\mathbf{y}(t+1,\omega),\cdots,\mathbf{y}(t+n,\omega)\right) \in X(t,\omega)$,

(i)  $\forall \omega \in \Omega, \ V\left(\mathbf{y}(t,\omega),\mathbf{y}(t+1,\omega),\cdots,\mathbf{y}(t+n,\omega),t,\omega\right) \in [-\infty,\infty)$,

(ii) the mapping $V\left(\mathbf{y}(t,\omega),\mathbf{y}(t+1,\omega),\cdots,\mathbf{y}(t+n,\omega),t,\omega\right): \Omega \to [-\infty,\infty)$ is measurable, and

(iii) $EV\left(\mathbf{y}(t,\omega),\mathbf{y}(t+1,\omega),\cdots,\mathbf{y}(t+n,\omega),t,\omega\right)$ exists in $[-\infty,\infty)$.

Suppose that the optimal path to (13) exists and is given by $\mathbf{y}^*(t,\omega)$, optimal in the sense of the overtaking criterion to be defined below.[10] We perturb it with a continuously differentiable curve $\mathbf{q}(t,\omega)$,

---

[9] Normally, $V$ is defined on $\left(\mathbb{R}^N\right)^{n+1} \times \mathbb{R} \times \Omega$. The domain of $V$ is denoted by $X(t,\omega)$, included in $\left(\mathbb{R}^N\right)^{n+1}$, for all $t$, $\omega$.

$$\mathbf{y}(t,\omega) = \mathbf{y}^*(t,\omega) + \varepsilon \cdot \mathbf{q}(t,\omega). \tag{14}$$

We define

$$\begin{aligned}&V(\varepsilon,T)\\&= \inf_{T \leq T'} \sum_{t=0}^{T'} [EV(\mathbf{y}^*(t,\omega) + \varepsilon \cdot \mathbf{q}(t,\omega), \mathbf{y}^*(t+1,\omega) + \varepsilon \cdot \mathbf{q}(t+1,\omega), \cdots, \mathbf{y}^*(t+n,\omega) + \varepsilon \cdot \mathbf{q}(t+n,\omega), t, \omega) \\&\quad - EV(\mathbf{y}^*(t,\omega), \mathbf{y}^*(t+1,\omega), \cdots, \mathbf{y}^*(t+n,\omega), t, \omega)].\end{aligned} \tag{15}$$

As in the preceding analysis, the optimal path satisfies the weak maximality criterion, á la Brock (1970), which is defined as:

**Definition 1'.** An attainable path $\{(\mathbf{y}^*(t,\omega)), t \in \mathbb{N}, \omega \in \Omega\}$ is optimal if no other attainable path overtakes it:

$$\lim_{T \to \infty} V(\varepsilon, T) \leq 0. \tag{16}$$

Let $V(\varepsilon) = \lim_{T \to \infty} V(\varepsilon, T)$. Differentiating it with respect to $\varepsilon$, we have

$$\begin{aligned}&\lim_{\varepsilon \to +0} \frac{V(\varepsilon)}{\varepsilon}\\&= \lim_{\varepsilon \to +0} \lim_{T \to \infty} \inf_{T \leq T'} \sum_{t=0}^{T'} \frac{1}{\varepsilon} [EV(\mathbf{y}^*(t,\omega) + \varepsilon \cdot \mathbf{q}(t,\omega), \cdots, \mathbf{y}^*(t+n,\omega) + \varepsilon \cdot \mathbf{q}(t+n,\omega), t, \omega) \\&\quad - EV(\mathbf{y}^*(t,\omega), \mathbf{y}^*(t+1,\omega), \cdots, \mathbf{y}^*(t+n,\omega), t, \omega)].\end{aligned} \tag{17}$$

Let $\lim_{\varepsilon \to +0} \frac{V(\varepsilon)}{\varepsilon} \equiv \Phi.$ We first assume

**Assumption 1'.** *Assume* $\Phi$ *converges uniformly for* $\varepsilon$ *when* $T \to \infty$.

---

[10] Again, our analysis only serves for the case where the optimal solution is interior to $X(t,\omega)$.

If $\Phi$ satisfies Assumption 1, we can then restate (17) as

$$\Phi = \lim_{T \to \infty} \liminf_{\varepsilon \to +0} \inf_{T \leq T'} \sum_{t=0}^{T'} \frac{1}{\varepsilon} [EV(\mathbf{y}^*(t,\omega) + \varepsilon \cdot \mathbf{q}(t,\omega), \cdots, \mathbf{y}^*(t+n,\omega) + \varepsilon \cdot \mathbf{q}(t+n,\omega), t, \omega) \quad (18)$$
$$- EV(\mathbf{y}^*(t,\omega), \mathbf{y}^*(t+1,\omega), \cdots, \mathbf{y}^*(t+n,\omega), t, \omega)].$$

We proceed to impose another assumption:

**Assumption 2'.** *We assume that for any $T > 0$,*

$$\inf_{T \leq T'} \sum_{t=0}^{T'} \frac{1}{\varepsilon} [EV(\mathbf{y}^*(t,\omega) + \varepsilon \cdot \mathbf{q}(t,\omega), \cdots, \mathbf{y}^*(t+n,\omega) + \varepsilon \cdot \mathbf{q}(t+n,\omega), t, \omega)$$
$$- EV(\mathbf{y}^*(t,\omega), \mathbf{y}^*(t+1,\omega), \cdots, \mathbf{y}^*(t+n,\omega), t, \omega)]$$

*converges uniformly for $\varepsilon$.*

Assume Assumptions 1' and 2', we can restate equality (18) as

$$\Phi = \liminf_{T \to \infty} \inf_{T \leq T'} \lim_{\varepsilon \to +0} \sum_{t=0}^{T'} \frac{1}{\varepsilon} [EV(\mathbf{y}^*(t,\omega) + \varepsilon \cdot \mathbf{q}(t,\omega), \cdots, \mathbf{y}^*(t+n,\omega) + \varepsilon \cdot \mathbf{q}(t+n,\omega), t, \omega)$$
$$- EV(\mathbf{y}^*(t,\omega), \mathbf{y}^*(t+1,\omega), \cdots, \mathbf{y}^*(t+n,\omega), t, \omega)].$$

(19)

Because $T'$ is finite uniformly for $\varepsilon$, if

$$\sum_{t=0}^{T'} \frac{1}{\varepsilon} [EV(\mathbf{y}^*(t,\omega) + \varepsilon \cdot \mathbf{q}(t,\omega), \cdots, \mathbf{y}^*(t+n,\omega) + \varepsilon \cdot \mathbf{q}(t+n,\omega), t, \omega)$$
$$- EV(\mathbf{y}^*(t,\omega), \mathbf{y}^*(t+1,\omega), \cdots, \mathbf{y}^*(t+n,\omega), t, \omega)]$$

exists, (19) is then rewritten as

$$\Phi = \liminf_{T \to \infty} \inf_{T \leq T'} \sum_{t=0}^{T'} \lim_{\varepsilon \to +0} \frac{1}{\varepsilon} [EV(\mathbf{y}^*(t,\omega) + \varepsilon \cdot \mathbf{q}(t,\omega), \cdots, \mathbf{y}^*(t+n,\omega) + \varepsilon \cdot \mathbf{q}(t+n,\omega), t, \omega)$$
$$- EV(\mathbf{y}^*(t,\omega), \mathbf{y}^*(t+1,\omega), \cdots, \mathbf{y}^*(t+n,\omega), t, \omega)].$$

(20)

We see that

$$\frac{1}{\varepsilon}[EV(\mathbf{y}^*(t,\omega)+\varepsilon\cdot\mathbf{q}(t,\omega),\cdots,\mathbf{y}^*(t+n,\omega)+\varepsilon\cdot\mathbf{q}(t+n,\omega),t,\omega)$$
$$-EV(\mathbf{y}^*(t,\omega),\mathbf{y}^*(t+1,\omega),\cdots,\mathbf{y}^*(t+n,\omega),t,\omega)]$$
$$=E[\frac{1}{\varepsilon}(V(\mathbf{y}^*(t,\omega)+\varepsilon\cdot\mathbf{q}(t,\omega),\cdots,\mathbf{y}^*(t+n,\omega)+\varepsilon\cdot\mathbf{q}(t+n,\omega),t,\omega)$$
$$-V(\mathbf{y}^*(t,\omega),\mathbf{y}^*(t+1,\omega),\cdots,\mathbf{y}^*(t+n,\omega),t,\omega))].$$

On the other hand, since $V$ is differentiable, we have

$$\lim_{\varepsilon\to+0}\left(\frac{V\left(\mathbf{y}^*(t,\omega)+\varepsilon\cdot\mathbf{q}(t,\omega),\cdots,\mathbf{y}^*(t+N-1,\omega)+\varepsilon\cdot\mathbf{q}(t+N-1,\omega),t,\omega\right)-V\left(\mathbf{y}^*(t,\omega),\cdots,\mathbf{y}^*(t+N-1,\omega),t,\omega\right)}{\varepsilon}\right)$$
$$=\sum_{i=1}^{n}[\frac{\partial V\left(\mathbf{y}^*(t,\omega),\mathbf{y}^*(t+1,\omega),\cdots,\mathbf{y}^*(t+N-1,\omega),t,\omega\right)}{\partial y_i(t,\omega)}q_i(t,\omega)$$
$$+\frac{\partial V\left(\mathbf{y}^*(t,\omega),\mathbf{y}^*(t+1,\omega),\cdots,\mathbf{y}^*(t+N-1,\omega),t,\omega\right)}{\partial y_i(t+1,\omega)}q_i(t+1,\omega)+\cdots$$
$$+\frac{\partial V\left(\mathbf{y}^*(t,\omega),\mathbf{y}^*(t+1,\omega),\cdots,\mathbf{y}^*(t+N-1,\omega),t,\omega\right)}{\partial y_i(t+N-1,\omega)}q_i(t+N-1,\omega)].$$

**Assumption 3'.** Define

$$m_t(\varepsilon,\omega)\equiv\frac{1}{\varepsilon}[V(\mathbf{y}^*(t,\omega)+\varepsilon\cdot\mathbf{q}(t,\omega),\cdots,\mathbf{y}^*(t+n,\omega)+\varepsilon\cdot\mathbf{q}(t+n,\omega),t,\omega)$$
$$-V(\mathbf{y}^*(t,\omega),\mathbf{y}^*(t+1,\omega),\cdots,\mathbf{y}^*(t+n,\omega),t,\omega)].$$

Assume that there exists $\hat{m}_t(\varepsilon,\omega)$ such that $|m_t(\tilde{\varepsilon},\omega)|<\hat{m}_t(\varepsilon,\omega)$, for any $\tilde{\varepsilon}\in(0,\varepsilon]$.

**Theorem 1'.** Under Assumption 1'-3', for any interior optimal path $\mathbf{y}^*(t,\omega)$, the Euler equation is

$$\frac{\partial V(0,\omega)}{\partial y_i(0,\omega)}=0,$$

$$\frac{\partial(V(0,\omega)+V(1,\omega))}{\partial y_i(1,\omega)}=0,$$

……

$$\frac{\partial(V(0,\omega)+\cdots+V(n-1,\omega))}{\partial y_i(t,\omega)}=0,$$

$$\frac{\partial \left(V(t-n,\omega)+\cdots+V(t,\omega)\right)}{\partial y_i(t,\omega)} = 0, \text{ for } n \leq t \leq T', \tag{21}$$

whereas the TVC is given by

$$\liminf_{T\to\infty} \sum_{T\leq T'} \sum_{i=1}^{n} \left( E\left( \frac{\partial \left(V(T'-n+1,\omega)+\cdots+V(T',\omega)\right)}{\partial y_i(T'+1,\omega)} q_i(T'+1,\omega) + \cdots + \frac{\partial \left(V(T',\omega)\right)}{\partial y_i(T'+n,\omega)} q_i(T'+n,\omega) \right) \right) \leq 0.$$

$$\tag{22}$$

**Proof.** Assume Assumption 3', from Lebesgue's dominated convergence theorem, we then have

$$\lim_{\varepsilon \to +0} E[\frac{1}{\varepsilon}(V(\mathbf{y}^*(t,\omega)+\varepsilon\cdot\mathbf{q}(t,\omega),\cdots,\mathbf{y}^*(t+n,\omega)+\varepsilon\cdot\mathbf{q}(t+n,\omega),t,\omega)$$
$$-V(\mathbf{y}^*(t,\omega),\mathbf{y}^*(t+1,\omega),\cdots,\mathbf{y}^*(t+n,\omega),t,\omega))]$$
$$= E\lim_{\varepsilon \to +0}[\frac{1}{\varepsilon}(V(\mathbf{y}^*(t,\omega)+\varepsilon\cdot\mathbf{q}(t,\omega),\cdots,\mathbf{y}^*(t+n,\omega)+\varepsilon\cdot\mathbf{q}(t+n,\omega),t,\omega)$$
$$-V(\mathbf{y}^*(t,\omega),\mathbf{y}^*(t+1,\omega),\cdots,\mathbf{y}^*(t+n,\omega),t,\omega))].$$

Hence,

$$\Phi = \liminf_{T\to\infty} \sum_{T\leq T'} \sum_{t=0}^{T'} \sum_{i=1}^{n} E[\frac{\partial V\left(\mathbf{y}^*(t,\omega),\mathbf{y}^*(t+1,\omega),\cdots,\mathbf{y}^*(t+n,\omega),t,\omega\right)}{\partial y_i(t,\omega)} q_i(t,\omega)$$
$$+ \frac{\partial V\left(\mathbf{y}^*(t,\omega),\mathbf{y}^*(t+1,\omega),\cdots,\mathbf{y}^*(t+n,\omega),t,\omega\right)}{\partial y_i(t+1,\omega)} q_i(t,\omega) + \cdots$$
$$+ \frac{\partial V\left(\mathbf{y}^*(t,\omega),\mathbf{y}^*(t+1,\omega),\cdots,\mathbf{y}^*(t+n,\omega),t,\omega\right)}{\partial y_i(t+n,\omega)} q_i(t+n,\omega)].$$

$$\tag{23}$$

We then obtain

$$\sum_{t=0}^{T'}\left(\sum_{i=1}^{n}E\left(\frac{\partial V(t,\omega)}{\partial y_i(t,\omega)}q_i(t,\omega)+\cdots+\frac{\partial V(t,\omega)}{\partial y_i(t+n,\omega)}q_i(t+n,\omega)\right)\right)$$

$$=E\{\sum_{i=1}^{n}\{\frac{\partial V(0,\omega)}{\partial y_i(0,\omega)}q_i(0,\omega)+\frac{\partial(V(0,\omega)+V(1,\omega))}{\partial y_i(1,\omega)}q_i(1,\omega)$$

$$+\cdots+\frac{\partial(V(0,\omega)+\cdots+V(n-1,\omega))}{\partial y_i(t,\omega)}q_i(n-1,\omega) \quad (24)$$

$$+\sum_{t=N-1}^{T'}\frac{\partial(V(t-n,\omega)+\cdots+V(t,\omega))}{\partial y_i(t,\omega)}q_i(t,\omega)$$

$$+\left(\frac{\partial(V(T'-n+1,\omega)+\cdots+V(T',\omega))}{\partial y_i(T'+1,\omega)}\right)q_i(T'+1,\omega)$$

$$+\cdots+\frac{\partial V(T',\omega)}{\partial y_i(T'+n,\omega)}q_i(T'+n,\omega)\}\}.$$

(21) and (22) can then be derived from (24).

*Q.E.D.*

Clearly, (22) is a stochastic version of the TVC in Cai and Nitta (2010). Again, when $\varepsilon \to {}^-0$, the argument is similar and the TVC is given as

$$\limsup_{T\to\infty}\sum_{T\le T'}^{n}\left(E\left(\frac{\partial(V(T'-n+1,\omega)+\cdots+V(T',\omega))}{\partial y_i(T'+1,\omega)}q_i(T'+1,\omega)+\cdots+\frac{\partial(V(T',\omega))}{\partial y_i(T'+n,\omega)}q_i(T'+n,\omega)\right)\right)\ge 0.$$

(22')

### 3. Examples and Counterexamples

#### 3.1 *An Example for the Continuous Time Problems*

We proceed to show that assumptions are imperative in the sense that (10) becomes invalid if one of them is violated. Assuming that $\Omega:\{1,\cdots,m\}$, $P:\Omega \to \mathbb{R}^+$, $\sum_{\omega=1}^{m} P(\omega)=1$, we consider a simple example:

$$v(x(t,\omega),\dot{x}(t,\omega),\ddot{x}(t,\omega),\omega) = (x(t,\omega)-\alpha(\omega))^2 + \beta \dot{x}(t,\omega) + \gamma \ddot{x}(t,\omega), \tag{27}$$

where $\alpha(\omega)>0$, $\beta(\omega)>0$, $\gamma(\omega)>0$ are constants, and the initial value $x_0$ is given, with $x_0 = \alpha(\omega)$. From (9), we see that the Euler equation is

$$v_1 - (v_2)' + (v_3)'' = 0, \tag{28}$$

that is,

$$2(x(t,\omega)-\alpha(\omega)) - (\beta(\omega))' + (\gamma(\omega))'' = 0. \tag{28'}$$

Thus, we have $x(t,\omega)=\alpha(\omega)$.

Choosing a $p(t,\omega)$ such that $p(0,\omega)=0$ and $p(t,\omega)>0$, $\dot{p}(0,\omega)=0$, there exists $T_0 > 0$ such that $\dot{p}(t,\omega)=0$, $t \geq T_0$, that is, $p(t,\omega)$ is a constant $p_\infty(\omega)>0$ for $t \geq T_0$. From (10), we see that the TVC is

$$\liminf_{T\to\infty} E\left[ p(t,\omega)\left(v_2(t,\omega)-(v_3(t,\omega))'\right) + \dot{p}(t,\omega)(v_3(t,\omega)) \right]_0^T \leq 0, \tag{29}$$

The left hand side (LHS) of (29) can be further rewritten as

$$\begin{aligned}
\text{LHS of (29)} &= \liminf_{T\to\infty} E\left[ p(t,\omega)(\beta(\omega)-\gamma'(\omega)) + \dot{p}(t,\omega)\gamma(\omega) \right]_0^T \\
&= \liminf_{T\to\infty} E\left((p(T,\omega)-p(0,\omega))(\beta(\omega)-\gamma'(\omega)) + \dot{p}(0,\omega)\gamma(\omega)\right) \\
&= \liminf_{T\to\infty} E\left(p_\infty(\beta(\omega)-\gamma'(\omega))\right) = E p_\infty(\omega)\beta(\omega) > 0.
\end{aligned}$$

We have then derived a contradiction to (10). Next we show that Assumption 1 is violated, which causes this contradiction.

We first consider

$$E(v(x(t,\omega)+\varepsilon p(t,\omega), \dot{x}(t,\omega)+\varepsilon \dot{p}(t,\omega), \ddot{x}(t,\omega)+\varepsilon \ddot{p}(t,\omega))-v(x(t,\omega), \dot{x}(t,\omega), \ddot{x}(t,\omega))).$$

Substituting $x(t,\omega)=\alpha(\omega)$ into it, we have

$$\begin{aligned}&E(v(x(t,\omega)+\varepsilon p(t,\omega), \dot{x}(t,\omega)+\varepsilon \dot{p}(t,\omega), \ddot{x}(t,\omega)+\varepsilon \ddot{p}(t,\omega))-v(x(t,\omega), \dot{x}(t,\omega), \ddot{x}(t,\omega)))\\ &=E((x(t,\omega)+\varepsilon p(t,\omega)-\alpha(\omega))^2+\beta(\omega)(\dot{x}(t,\omega)+\varepsilon \dot{p}(t,\omega))+\gamma(\omega)(\ddot{x}(t,\omega)+\varepsilon \ddot{p}(t,\omega))\\ &\quad -((x(t,\omega)-\alpha(\omega))^2+\beta(\omega)(\dot{x}(t,\omega))+\gamma(\omega)(\ddot{x}(t,\omega))))\\ &=E((\varepsilon p(t,\omega))^2+\varepsilon(\beta(\omega)\dot{p}(t,\omega)+\gamma(\omega)\ddot{p}(t,\omega))).\end{aligned}$$

(30)

Hence,

$$\begin{aligned}&\inf_{T\le T'}\int_0^{T'}E\left(\frac{(\varepsilon p(t,\omega))^2+\varepsilon(\beta(\omega)\dot{p}(t,\omega)+\gamma(\omega)\ddot{p}(t,\omega))}{\varepsilon}\right)dt\\ &=\inf_{T\le T'}\int_0^{T'}E(\varepsilon p(t,\omega)^2+(\beta(\omega)\dot{p}+\gamma(\omega)\ddot{p}))dt\\ &=\inf_{T\le T'}E\left(\varepsilon\int_0^{T'}p(t,\omega)^2\,dt+[\beta(\omega)p(t,\omega)+\gamma(\omega)\dot{p}(t,\omega)]_0^{T'}\right)\\ &=\inf_{T\le T'}\left(\varepsilon\int_0^{T'}Ep(t,\omega)^2\,dt+E(\beta(\omega)(p(T',\omega)-p(0,\omega))+\gamma(\omega)(\dot{p}(T',\omega)-\dot{p}(0,\omega)))\right)\\ &=\inf_{T\le T'}\left(\varepsilon\int_0^{T'}Ep(t,\omega)^2\,dt+E\beta(\omega)p_\infty(\omega)\right)\\ &=\inf_{T\le T'}\left(\varepsilon\int_0^{T'}Ep(t,\omega)^2\,dt\right)+bEp_\infty(\omega)\\ &=\varepsilon\int_0^{T}Ep(t,\omega)^2\,dt+\beta(\omega)Ep_\infty(\omega).\end{aligned}\qquad(31)$$

$\Pi$ is the limit of (31) when $T\to\infty$, $\varepsilon\to 0$. We first assume that $\int_0^\infty p(t,\omega)^2\,dt=\infty.$ Because $\lim_{\varepsilon\to\infty}\lim_{T\to\infty}\varepsilon\int_0^T Ep(t,\omega)^2\,dt=\infty,$ whereas $\lim_{T\to\infty}\lim_{\varepsilon\to 0}\varepsilon\int_0^T Ep(t,\omega)^2\,dt=0$, we see that $\Pi$ does not converge uniformly for $\varepsilon$ when $T \to \infty$. Hence, Assumption 1 is violated and (10) (when $n = 2$, $\liminf\limits_{T \to \infty} E\left(x^*(t,\omega)(v_2(t,\omega) - (v_3(t,\omega))') + (\dot{x}^*(t,\omega))v_3(t,\omega)\right)\big|_T \leq 0$) is also not satisfied.

Alternatively, we assume that $\int_0^\infty p(t,\omega)^2 dt < \infty$. Because now $\lim\limits_{\varepsilon \to \infty} \lim\limits_{T \to \infty} \varepsilon \int_0^T Ep(t,\omega)^2 dt = \lim\limits_{T \to \infty} \lim\limits_{\varepsilon \to 0} \varepsilon \int_0^T Ep(t,\omega)^2 dt$, we see that $\Pi$ converges uniformly for $\varepsilon$ when $T \to \infty$. Under such a case, both Assumption 1 and (10) are satisfied.

In sum, we see that even for the same unbounded $U$, depending on the properties of the perturbing curve $p(t,\omega)$, Assumption 1 can be either satisfied or violated. This also holds for Assumption 2. In other words, depending on the properties of the perturbation curves under consideration, the example can also represent a counterexample, in which the assumptions are violated and the TVC is not satisfied.

### 3.2  An Example for the Discrete Time Problems

We assume that $\Omega:\{1,\cdots,m\}$, $P:\Omega \to \mathbb{R}^+$, $\sum\limits_{\omega=1}^m P(\omega) = 1$, and consider the following example:

$$V(t,\omega) = V(y(t,\omega), y(t+1,\omega), y(t+2,\omega), t, \omega) \\ = (y(t,\omega) - \alpha(\omega))^2 + \beta(\omega) y(t+1,\omega) + \gamma(\omega) y(t+2,\omega), \tag{32}$$

where $\alpha(\omega) > 0$, $\beta(\omega) > 0$, $\gamma(\omega) > 0$, and the initial value $y(0,\omega) = y_0(\omega)$, $y(1,\omega) = y_1(\omega)$ is given. From (21), we see that Euler equation is given by

$$\frac{\partial V(0,\omega)}{\partial y(0,\omega)} = 0,$$

$$\frac{\partial (V(0,\omega) + V(1,\omega))}{\partial y(1,\omega)} = 0,$$

$$\frac{\partial (V(t-2,\omega) + V(t-1,\omega) + V(t,\omega))}{\partial y(t,\omega)} = 0, \quad 2 \leq t \leq T', \tag{33}$$

which implies

$$t = 0, \; 2(y(0,\omega) - \alpha(\omega)) = 0$$

$$t = 1, \; 2(y(1,\omega) - \alpha(\omega)) + \beta(\omega) = 0$$

$$t = 2, \; 2(y(2,\omega) - \alpha(\omega)) + \beta(\omega) + \gamma(\omega) = 0,$$

$$t = 3, \; 2(y(3,\omega) - \alpha(\omega)) + \beta(\omega) + \gamma(\omega) = 0,$$

……,

$$t = T', \; 2(y(T',\omega) - \alpha(\omega)) + \beta(\omega) + \gamma(\omega) = 0. \tag{33'}$$

Thus, we have

$$y(0,\omega) = \alpha(\omega),$$

$$y(1,\omega) = \alpha(\omega) - \frac{\beta(\omega)}{2},$$

$$y(2,\omega) = y(3,\omega) = \cdots = y(T',\omega) = \alpha(\omega) - \frac{\beta(\omega) + \gamma(\omega)}{2}.$$

Choosing a $q(t,\omega)$ so that $q(0,\omega) = 0$ and $q(t,\omega) > 0$, there exists $T_0 > 0$, $q(t,\omega)$ is a constant $q_\infty(\omega) > 0$ when $t \geq T_0$.

From (32), we see that

$$\frac{\partial (V(T'-1,\omega) + V(T',\omega))}{\partial y(t+1,\omega)} q(T'+1,\omega) = (\gamma(\omega) + \beta(\omega)) q(T'+1,\omega) > 0, \tag{34}$$

$$\frac{\partial (V(T',\omega))}{\partial y(t+2,\omega)} q(T'+2,\omega) = \gamma(\omega) q(T'+2,\omega) > 0. \tag{35}$$

Hence, we derived a contradiction to (22).

Next we show that Assumption 1' is violated, which causes this contradiction. We consider

$$EV(y(t,\omega)+\varepsilon q(t,\omega), y(t+1,\omega)+\varepsilon q(t+1,\omega), y(t+2,\omega)+\varepsilon q(t+2,\omega),t,\omega)$$
$$-EV(y(t,\omega), y(t+1,\omega), y(t+2,\omega),t,\omega).$$

Substituting $y(t,\omega) = \alpha(\omega) - \dfrac{\beta(\omega)+\gamma(\omega)}{2}$ into it, we have

$$EV(y(t,\omega)+\varepsilon q(t,\omega), y(t+1,\omega)+\varepsilon q(t+1,\omega), y(t+2,\omega)+\varepsilon q(t+2,\omega),t,\omega)$$
$$-EV(y(t,\omega), y(t+1,\omega), y(t+2,\omega),t,\omega)$$
$$= E(y(t,\omega)+\varepsilon q(t,\omega)-\alpha(\omega))^2 + \beta(\omega)(y(t+1,\omega)+\varepsilon q(t+1,\omega)) + \gamma(\omega)(y(t+2,\omega)+\varepsilon q(t+2,\omega))$$
$$- E\left((y(t,\omega)-\alpha(\omega))^2 + \beta(\omega)y(t+1,\omega) + \gamma(\omega)y(t+2,\omega)\right)$$
$$= E\left(\varepsilon q(t,\omega) - \left(\dfrac{\beta(\omega)+\gamma(\omega)}{2}\right)\right)^2 + \varepsilon(\beta(\omega)q(t+1,\omega)+\gamma(\omega)q(t+2,\omega)) - \left(\dfrac{\beta(\omega)+\gamma(\omega)}{2}\right)^2.$$

(36)

Hence,

$$\inf_{T \leq T'} \sum_{t=0}^{T'} E\left(\dfrac{1}{\varepsilon}\left(\varepsilon q(t,\omega) - \left(\dfrac{\beta(\omega)+\gamma(\omega)}{2}\right)\right)^2 + \varepsilon\beta(\omega)q(t+1,\omega) + \varepsilon\gamma(\omega)q(t+2,\omega) - \left(\dfrac{\beta(\omega)+\gamma(\omega)}{2}\right)^2\right)$$

$$= \inf_{T \leq T'} \sum_{t=0}^{T'} E\left(\varepsilon q(t,\omega)^2 - (\beta(\omega)+\gamma(\omega))q(t,\omega) + \beta(\omega)q(t+1,\omega) + \gamma(\omega)q(t+2,\omega)\right)$$

$$= \inf_{T \leq T'} E\left(\varepsilon \sum_{t=0}^{T'} q(t,\omega)^2 + \beta(\omega)q(T'+1,\omega) + \gamma(\omega)q(T'+2,\omega)\right)$$

$$= \inf_{T \leq T'} E\left(\left(\varepsilon \sum_{t=0}^{T'}(q(t,\omega)^2)\right) + \beta(\omega)q_\infty(\omega) + \gamma(\omega)q_\infty(\omega)\right)$$

$$= E\varepsilon \sum_{t=0}^{T}(q(t,\omega)^2) + \beta(\omega)q_\infty(\omega) + \gamma(\omega)q_\infty(\omega).$$

(37)

$\Phi$ is the limit of (37) when $T \to \infty$, $\varepsilon \to 0$. Similar to the example in 3.1, when

$$\sum_{t=0}^{\infty}(q(t,\omega)^2) = \infty, \quad \text{because} \quad \lim_{\varepsilon \to 0}\lim_{T \to \infty} E\left(\varepsilon \sum_{t=0}^{T}(q(t,\omega)^2) + \beta(\omega)q_\infty(\omega) + \gamma(\omega)q_\infty(\omega)\right) = \infty ,$$

whereas

$$\lim_{T\to\infty}\lim_{\varepsilon\to 0} E\left(\varepsilon\sum_{t=0}^{T}\left(q(t,\omega)^2\right)+\beta(\omega)q_\infty(\omega)+\gamma(\omega)q_\infty(\omega)\right)=E\left(\beta(\omega)q_\infty(\omega)+\gamma(\omega)q_\infty(\omega)\right),$$

we see that $\Phi$ does not converge uniformly for $\varepsilon$ when $T\to\infty$. Hence, Assumption 1' is violated and (22) is also not satisfied. On the other hand, when $\sum_{t=0}^{\infty}\left(q(t,\omega)^2\right)<\infty$, because

$$\lim_{\varepsilon\to 0}\lim_{T\to\infty} E\left(\varepsilon\sum_{t=0}^{T}\left(q(t,\omega)^2\right)+\beta(\omega)q_\infty(\omega)+\gamma(\omega)q_\infty(\omega)\right)$$
$$=\lim_{T\to\infty}\lim_{\varepsilon\to 0} E\left(\varepsilon\sum_{t=0}^{T}\left(q(t,\omega)^2\right)+\beta(\omega)q_\infty(\omega)+\gamma(\omega)q_\infty(\omega)\right),$$

we see that $\Phi$ converges uniformly for $\varepsilon$ when $T\to\infty$. Under such a case, both Assumption 1' and (22) are satisfied. Hence, similar to the differential case, depending on the properties of the perturbation curves under consideration, the example also represents a counterexample.

## 4. Discussion

In this section, we investigate the correspondence between the results for the continuous time models and those for the discrete time models. We only consider the third-order case.

For the first order case, because $\dot{x}(t)$ corresponds to $(x_{t+1}-x_t)$, $v(x,\dot{x},\ddot{x},t,\omega)$ would correspond to $V(y(t,\omega),y(t+1,\omega)-y(t,\omega),y(t+2,\omega)-y(t+1,\omega),t,\omega)$. We then define

$$v(x,y,z,t,\omega)=V(x,x+y,x+2y+z,t,\omega).$$

Clearly, $v_1 = V_1 + V_2 + V_3$, $v_2 = V_2 + 2V_3$, $v_3 = V_3$.

Substituting the above into (21), we see that the Euler equation can be stated as

$$\frac{\partial(V(t,\omega)+V(t+1,\omega)+V(t+2,\omega))}{\partial c(t+2,\omega)}$$
$$= V_3(t,\omega) + V_2(t+1,\omega) + V_1(t+2,\omega)$$
$$= v_3(t,\omega) + v_2(t+1,\omega) - 2V_3(t+1,\omega) + v_1(t+2,\omega) - V_2(t+2,\omega) - V_3(t+2,\omega)$$
$$= v_3(t,\omega) + v_2(t+1,\omega) - 2v_3(t+1,\omega) + v_1(t+2,\omega) - v_2(t+2,\omega) + 2v_3(t+2,\omega) - v_3(t+2,\omega)$$
$$= v_3(t,\omega) + v_2(t+1,\omega) - 2v_3(t+1,\omega) + v_1(t+2,\omega) - v_2(t+2,\omega) + v_3(t+2,\omega)$$
$$= v_1(t+2,\omega) + (v_2(t+1,\omega) - v_2(t+2,\omega)) + (v_3(t,\omega) - 2v_3(t+1,\omega) + v_3(t+2,\omega))$$
$$= 0,$$

which clearly corresponds to $v_1 - \dot{v}_2 + \ddot{v}_3 = 0$.

Next we consider the correspondence between TVCs. The discrete time third-order TVC is given by

$$\frac{\partial(V(T'-1,\omega)+V(T',\omega))}{\partial c(T'+1,\omega)} q(T'+1,\omega) + \frac{\partial(V(T',\omega))}{\partial c(T'+2,\omega)} q(T'+2,\omega)$$
$$= (V_3(T'-1,\omega) + V_2(T',\omega)) q(T'+1,\omega) + V_3(T',\omega) q(T'+2,\omega)$$
$$= (v_3(T'-1,\omega) + v_2(T',\omega) - 2v_3(T',\omega)) q(T'+1,\omega) + v_3(T,\omega) q(T'+2,\omega)$$
$$= (v_3(T'-1,\omega) - (v_3(T',\omega) - v_2(T',\omega))) q(T'+1,\omega) - v_3(t,\omega) q(T'+1,\omega) + v_3(T,\omega) q(T'+2,\omega)$$
$$= q(T'+1,\omega)(v_3(T'-1,\omega) - (v_3(T',\omega) - v_2(T',\omega))) + (q(T'+2,\omega) - q(T'+1,\omega)) v_3(T,\omega)$$
$$= 0,$$

which corresponds to $p\left(v_2 - (v_3)'\right) + \dot{p}(v_3) = 0$.

Hence, the Euler equations and the TVCs for the discrete time models correspond to those for the continuous time models. It has been shown that, for the first-order case, the results for the Euler equations and TVCs are intrinsically the same for the continuous time models and discrete time models (Kamihigashi, 2004). Our results suggest that this property also holds for stochastic higher-order problems.

# 5. Application: Simplified Household Maximization Problem in Futagami and Iwaisako (2007)

Consider the following simplified household maximization problem in Futagami and Iwaisaka (2007):

$$\max \sum_{t=0}^{\infty} V(t),$$

where $V(t) \equiv V(y_t, y_{t+1}, \cdots, y_{t+n}, t) = \begin{cases} 0, t \leq n-1 \\ \beta^t \ln(c_t), t \geq n \end{cases}$, $c_t = -y_{t+n} + y_{t+n-1} + \cdots + y_{t+1} + y_t$.

Clearly, $V(0) = 0$, $V(1) = 0$, and $V(n-1) = 0$. From Theorem 1', the Euler equation is given by

$$\frac{\partial (V(t-n) + \cdots + V(t))}{\partial y_t} = 0, \quad n \leq t \leq T,$$

which can then be written as

$$\frac{1}{c_{t-n}} = \beta \frac{1}{c_{t-n+1}} + \beta^2 \frac{1}{c_{t-n+2}} + \cdots + \beta^n \frac{1}{c_t}.$$

On the other hand, the TVC is given by

$$\liminf_{T \to \infty, T \leq T'} \left( \left( \frac{1}{c_{T'-n+2}} + \beta^1 \frac{1}{c_{T'-n+2}} + \cdots + \beta^{n-2} \frac{1}{c_{T'}} \right) q_{T'+1} + \cdots + \beta^{n-2} \frac{1}{c_{T'}} q_{T'+n-1} \right) = 0.$$

Clearly, for arbitrary $U(c_t, t)$, where there exists $\frac{\partial U(c_t, t)}{\partial c_t}$ and $V(y_t, \cdots, y_{t+n}, t) \equiv U(c_t, t)$, we can easily extend the above analysis under Assumptions 1', 2', 3' for $U(c_t, t)$, $y_t^*$ and $q_t$.

# 6. Concluding Remarks

In this paper, we present three assumptions that would be needed, in addition to the standard assumptions, when examining stochastic infinite horizon optimization problems with unbounded objective functions. Our results generalize the results of Ekeland and Scheinkman (1986), Michel (1990), and Kamihigashi (2001, 2003), which considered first-order problems, to stochastic higher-order problems. Moreover, our results also extend the results of Okumura, Cai and Nitta (2009) and Cai and Nitta (2010), which considered higher-order problems, to stochastic cases. Specifically, Assumptions 1 and 2, and Assumptions 1' and 2' constitute the stochastic versions of Assumptions 1 and 2 presented in Okumura, Cai and Nitta (2009) and Cai and Nitta (2010), respectively, whereas Assumption 3' corresponds to Assumption (A.4) in Kamihigashi (2003), which considered stochastic first-order difference problems. Clearly, our assumptions hold when a discounting factor is incorporated into the model. In this sense, this paper also generalizes the TVCs examined in the presence of discounting.